\documentclass[12pt]{article}
\usepackage{graphicx}
\usepackage{amsmath,amsthm,amssymb,enumerate}%, esint}
\usepackage{euscript,mathrsfs}
\usepackage[left=2cm,right=2cm,top=3.5cm,bottom=3.5cm]{geometry}
\usepackage{color}
\catcode`\@=11 \@addtoreset{equation}{section}

\catcode`\@=12

\allowdisplaybreaks

\newtheorem{Theorem}{Theorem}[section]
\newtheorem{Proposition}[Theorem]{Proposition}
\newtheorem{Lemma}[Theorem]{Lemma}
\newtheorem{Corollary}[Theorem]{Corollary}

\theoremstyle{definition}
\newtheorem{Definition}[Theorem]{Definition}

\newtheorem{Remark}[Theorem]{Remark}

\newcommand{\bTheorem}[1]{
\begin{Theorem} \label{T#1} }
\newcommand{\eT}{\end{Theorem}}

\newcommand{\bProposition}[1]{
\begin{Proposition} \label{P#1}}
\newcommand{\eP}{\end{Proposition}}

\newcommand{\bLemma}[1]{
\begin{Lemma} \label{L#1} }
\newcommand{\eL}{\end{Lemma}}

\newcommand{\bCorollary}[1]{
\begin{Corollary} \label{C#1} }
\newcommand{\eC}{\end{Corollary}}

\newcommand{\bRemark}[1]{
\begin{Remark} \label{R#1} }
\newcommand{\eR}{\end{Remark}}

\newcommand{\bDefinition}[1]{
\begin{Definition} \label{D#1} }
\newcommand{\eD}{\end{Definition}}

\newcommand{\vcg}[1]{{\pmb #1}}

\newcommand{\bFormula}[1]{
\begin{equation} \label{#1}}
\newcommand{\eF}{\end{equation}}

\newcommand{\Ov}[1]{\overline{#1}}

\newcommand{\DC}{C^\infty_c}

\newcommand{\aleq}{\stackrel{<}{\sim}}
\newcommand{\ageq}{\stackrel{>}{\sim}}

\newcommand{\vr}{\varrho}
\newcommand{\vre}{\vr_\ep}

\newcommand{\vte}{\vt_\ep}
\newcommand{\vue}{\vu_\ep}

\newcommand{\vt}{\vartheta}
\newcommand{\vu}{\vc{u}}
\newcommand{\vc}[1]{{\bf #1}}

\newcommand{\Div}{{\rm div}_x}
\newcommand{\Grad}{\nabla_x}
\newcommand{\Dt}{\frac{\rm d}{{\rm d}t}}

\newcommand{\dx}{\,{\rm d} {x}}

\newcommand{\dt}{\,{\rm d} t }

\newcommand{\intO}[1]{\int_{\Omega} #1 \ \dx}

\newcommand{\ep}{\varepsilon}

\definecolor{Cgrey}{rgb}{0.85,0.85,0.85}
\definecolor{Cblue}{rgb}{0.50,0.85,0.85}
\definecolor{Cred}{rgb}{1,0,0}
\definecolor{fancy}{rgb}{0.10,0.85,0.10}

\newcommand\Cbox[2]{%
    \newbox\contentbox%
    \newbox\bkgdbox%
    \setbox\contentbox\hbox to \hsize{%
        \vtop{
            \kern\columnsep
            \hbox to \hsize{%
                \kern\columnsep%
                \advance\hsize by -2\columnsep%
                \setlength{\textwidth}{\hsize}%
                \vbox{
                    \parskip=\baselineskip
                    \parindent=0bp
                    #2
                }%
                \kern\columnsep%
            }%
            \kern\columnsep%
        }%
    }%
    \setbox\bkgdbox\vbox{
        \color{#1}
        \hrule width  \wd\contentbox %
               height \ht\contentbox %
               depth  \dp\contentbox
        \color{black}
    }%
    \wd\bkgdbox=0bp%
    \vbox{\hbox to \hsize{\box\bkgdbox\box\contentbox}}%
    \vskip\baselineskip%
}

\date{}

\begin{document}

%%%%%%%%%%%%%%%%%%%%%%%%%%%%%%%%

\title{Measure-valued solutions to the complete Euler system}

\author{Jan B\v rezina \and Eduard Feireisl
\thanks{The research of E.F.~leading to these results has received funding from the
European Research Council under the European Union's Seventh
Framework Programme (FP7/2007-2013)/ ERC Grant Agreement
320078. The Institute of Mathematics of the Academy of Sciences of
the Czech Republic is supported by RVO:67985840.}
}

\date{\today}

\maketitle

\bigskip

\centerline{Tokyo Institute of Technology}

\centerline{ 2-12-1 Ookayama, Meguro-ku, Tokyo, 152-8550, Japan}

\bigskip

\centerline{Institute of Mathematics of the Academy of Sciences of the Czech Republic}

\centerline{\v Zitn\' a 25, CZ-115 67 Praha 1, Czech Republic}

\bigskip

\begin{abstract}

We introduce the concept of \emph{dissipative measure-valued solution} to the complete Euler system describing the motion of an inviscid compressible fluid.
These solutions are characterized by a parameterized (Young) measure and a dissipation defect in the total energy balance. The dissipation defect dominates the concentration
errors in the equations satisfied by the Young measure. A dissipative measure-valued solution can be seen as the most general concept of solution to the Euler
system retaining its structural stability. In particular, we show that a dissipative measure-valued solution necessarily coincides with a classical one
on its life span provided they share the same initial data.

\end{abstract}


%\tableofcontents

\section{Introduction}
\label{i}

In his pioneering work \cite{DiP2}, DiPerna proposed a new concept of solution, known as \emph{measure--valued solution}, to nonlinear systems of partial differential equations
admitting uncontrollable oscillations. In particular with focus on the compressible Euler system and other related models of \emph{inviscid} fluids. Later on,
a similar strategy has been adopted even to problems involving \emph{viscous fluid flows}, where compactness of the solution set is either absent or out of reach of
the available mathematical tools, see e.g. the monograph Ne\v cas et al. \cite{MNRR} and the references therein.
Although \emph{existence} of a measure-valued solution to a given problem is usually an almost straightforward consequence of {\it a priori} bounds, its \emph{uniqueness} in terms of the initial data can be seen as the weakest point of this approach. In addition, the
recent results of DeLellis, Sz\' ekelyhidi and their collaborators \cite{Chiod}, \cite{ChiDelKre}, \cite{DelSze3} show that {uniqueness} may be in fact
violated even within the class of more conventional weak solutions satisfying the standard entropy admissibility criteria.

Brenier et al. \cite{BrDeSz} proposed a new approach seeing the measure-valued solutions as possibly the largest class in which the family of  smooth (classical) solutions
is stable. In particular, they show the so-called weak (measure-valued)--strong uniqueness principle
for the incompressible Euler system. Specifically,
a classical and a measure--valued solution emanating from the same initial data coincide as long as the former exists. These results have been extended
to the isentropic Euler and Navier--Stokes systems by Gwiazda et al. \cite{GSWW} and \cite{FGSWW1}. The recently renewed interest in measure--valued
solutions in fluid mechanics has been also initiated by certain numerical experiments with oscillatory solutions, see Fjordholm et al. \cite{FjKaMiTa}, \cite{FjMiTa2}, \cite{FjMiTa1}. Following the philosophy of Brenier et al. \cite{BrDeSz}, we focus on the concept of measure-valued solutions in the widest possible sense. Accordingly, using the fundamental laws of
thermodynamics, we extract the minimal piece of information to be retained to preserve the weak-strong uniqueness principle.

We consider the \emph{complete Euler system} describing the time evolution of the mass density $\vr = \vr(t,x)$, the velocity $\vu = \vu(t,x)$, and the (absolute)
temperature $\vt = \vt(t,x)$ of a compressible inviscid fluid:

\Cbox{Cgrey}{

\begin{eqnarray}
\label{i1} \partial_t \vr + \Div (\vr \vu) &=& 0, \\
\label{i2} \partial_t (\vr \vu) + \Div (\vr \vu \otimes \vu) + \Grad p(\vr, \vt) &=&0, \\
\label{i3} \partial_t \left( \frac{1}{2} \vr |\vu|^2 + \vr e( \vr, \vt) \right) +
\Div \left[ \left( \frac{1}{2} \vr |\vu|^2 + \vr e( \vr, \vt) + p(\vr, \vt) \right) \vu \right] &=& 0,\\
\label{i4} \partial_t (\vr s(\vr, \vt)) + \Div (\vr s(\vr,\vt) \vu ) &\geq& 0,
\end{eqnarray}

}

\noindent
where the pressure $p = p(\vr,\vt)$, the specific internal energy $e = e(\vr, \vt)$, and the specific entropy $s = s(\vr, \vt)$ are interrelated through \emph{Gibbs' equation}
\begin{equation} \label{i5}
\vt D s(\vr, \vt) = D e(\vr, \vt) + p(\vr, \vt) D \left( \frac{1}{\vr} \right).
\end{equation}

If $p,e,s$ comply with (\ref{i5}), then any \emph{smooth} solution of (\ref{i1}--\ref{i3}) satisfies automatically the entropy balance
in (\ref{i4}),
\begin{equation} \label{i6}
\partial_t (\vr s(\vr, \vt)) + \Div (\vr s(\vr,\vt) \vu ) = 0.
\end{equation}
This is no longer true for
the weak solutions, here typically represented by shock waves, for which the entropy inequality (\ref{i4}) may be appended to the
weak formulation of (\ref{i1}--\ref{i3}) as an admissibility criterion imposed by the Second law of thermodynamics, see e.g. the monograph
Benzoni-Gavage and Serre \cite{BenSer}.

Our goal is to address the problem of weak (measure-valued)--strong uniqueness for the Euler system (\ref{i1}--\ref{i4}). Accordingly, we focus
on identifying the \emph{largest class possible} of measure--valued solutions, in which such a result holds, rather than the optimal one with respect to the expected regularity  of solutions. To this end, we follow the approach advocated in \cite{FENO6}, where equations (\ref{i1}), (\ref{i2}), with \emph{inequality} (\ref{i4}),  are
supplemented with the total energy inequality
\begin{equation} \label{i7}
\Dt \intO{ \left[ \frac{1}{2} \vr |\vu|^2 + \vr e( \vr, \vt) \right] } \leq 0,
\end{equation}
where $\Omega \subset R^3$ is the physical domain occupied by the fluid. To simplify presentation, we impose the periodic boundary conditions, meaning $\Omega$ can be identified with the flat torus
\[
\Omega = \left( [0,1] |_{\{ 0,1 \}} \right)^3.
\]
The problem is closed by prescribing the initial data
\begin{equation} \label{i8}
\vr(0, \cdot) = \vr_0, \ \vt(0, \cdot) = \vt_0, \ \vu(0, \cdot) = \vu_0.
\end{equation}

\begin{Remark} \label{rr11}

We focus on the most difficult and physically relevant case of the spatial dimension $N = 3$. The arguments can be easily adapted to
$N=1,2$ as well.

\end{Remark}

The paper is organized as follows. In Section \ref{MV}, we define the measure valued solutions motivated by considering the cluster points of (hypothetical)
families of weak solutions in the topology induced by the available {\it a priori} bounds. We point out that our class of measure-valued solutions is larger
than that one proposed by Kr\" oner and Zajaczkowski \cite{KrZa}, where the entropy \emph{equality} (\ref{i6}) is required. In particular, it includes all
admissible weak solutions to the Euler system. Another conceptually new feature of our approach is that the \emph{dissipation defect}
``hidden'' behind the inequality sign in (\ref{i4}) and (\ref{i7}) dominates the concentration error emerging in (\ref{i2}).

In Section \ref{WS}, we establish and prove the main result of the present paper - the weak (measure-valued)--strong uniqueness principle.
To this end, we use the relative energy inequality for system (\ref{i1}, \ref{i2},
\ref{i4}, \ref{i7}) identified in \cite{FeiNov10}. One of the principal difficulties is the hypothetical presence of vacuum zones on which
the underlying equations fail to provide any control on the behavior of solutions. To avoid this problem, new phase variables must be considered - the density $\vr$,
the internal energy density $E = \vr e$, and the momentum $\vc{m} = \vr \vu$.

Finally, possible extensions and applications of the main result are discussed in Section \ref{C}.

\section{Measure-valued solutions}
\label{MV}

Motivated by \cite{FGSWW1}, \cite{GSWW}  we introduce  the concept of \emph{dissipative measure--valued  solution} to the Euler system. For the sake of simplicity, we start with the constitutive equations of a perfect gas, specifically
\begin{equation} \label{MV1}
p(\vr, \vt) =  \vr \vt, \ e(\vr, \vt) = c_v \vt, \  s(\vr, \vt) = \log\left( \frac{\vt^{c_v}}{\vr} \right),
\end{equation}
where $c_v > 0$ is the (constant) specific heat at constant volume.

In the seminal work of DiPerna \cite{DiP2}, the measure-valued solutions have been identified as (weak) limits of weak solutions of system (\ref{i1}--\ref{i4}) or its suitable viscous approximation. In particular, all weak solutions of the problem should fall into the category of measure-valued solutions.

A weak solution $[\vr, \vt, \vu]$ of the Euler system
in $(0,T) \times \Omega$, supplemented with the initial data (\ref{i8}), satisfies the family of integral identities:

\begin{equation} \label{MV2}
\int_0^T \intO{ \left[ \vr \partial_t \varphi + \vr \vu \cdot \Grad \varphi \right] } \dt = - \intO{ \vr_0 \varphi(0, \cdot) },
\end{equation}
for any $\varphi \in \DC([0,T) \times \Omega)$;
\begin{equation} \label{MV3}
\int_0^T \intO{ \left[ \vr \vu \cdot \partial_t\vcg{\varphi} + \vr \vu \otimes \vu : \Grad \vcg{\varphi} + \vr \vt \Div \vcg{\varphi} \right] } \dt =
- \intO{ \vr_0 \vu_0 \cdot \vcg{\varphi}(0, \cdot) },
\end{equation}
for any $\vcg{\varphi} \in \DC([0,T) \times \Omega; R^3)$;
\begin{equation} \label{MV4}
\begin{split}
\int_0^T \intO{ \left[ \left( \frac{1}{2} \vr |\vu|^2 + c_v \vr \vt \right) \partial_t \varphi +
 \left( \frac{1}{2} \vr |\vu|^2 + c_v \vr \vt \right) \vu \cdot \Grad \varphi + \vr \vt \vu \cdot \Grad \varphi \right] } \dt \\
= - \intO{ \left( \frac{1}{2} \vr_0 |\vu_0|^2 + c_v \vr_0 \vt_0 \right) \varphi(0, \cdot) },
\end{split}
\end{equation}
for any $\varphi \in \DC([0,T) \times \Omega)$.
In addition, a weak solution is called \emph{admissible} if the entropy inequality is imposed
\begin{equation} \label{MV5}
\begin{split}
\int_0^T \intO{ \left[ \vr Z\left( \log \left( \frac{ \vt^{c_v} }{\vr} \right) \right)  \partial_t \varphi + \vr Z \left( \log \left( \frac{ \vt^{c_v} }{\vr} \right) \right) \vu \cdot \Grad \varphi  \right] } \dt \\ \leq -
\intO{ \vr_0 Z \left( \log \left( \frac{ \vt_0^{c_v} }{\vr_0} \right) \right) \varphi (0, \cdot) },
\end{split}
\end{equation}
for any $\varphi \in \DC([0,T) \times \Omega)$, $\varphi \geq 0$, $Z \in BC(R)$, $Z' \geq 0$.

The weak solutions should also satisfy the natural constraints
\[
\vr(t,x) \geq 0, \ \vt (t,x) > 0 \ \mbox{a.a. in} \ (0,T) \times \Omega,
\]
in particular, the entropy density $\vr s(\vr, \vt)$ is well defined.

\begin{Remark} \label{bfr1}
The use of the cut-off function $Z$ in (\ref{MV5}) is motivated by Chen and Frid \cite{CheFr2}. Inequality (\ref{MV5}) may be seen as a renormalized
version of (\ref{i4}). 

\end{Remark}

As is well known, smooth solutions of the Euler system may develop singularities in finite time for a fairly general class of initial data. The admissible weak solutions
represent a physically grounded alternative providing the description of the system in an arbitrary time lap. Unfortunately, global-in-time existence of weak solutions
is a largely open problem. In addition, the recent examples provided by the theory of convex integration, see Chiodaroli et al. \cite{Chiod}, \cite{ChiDelKre}, show that
uniqueness may fail in the multidimensional case even in the class of admissible weak solutions.

\subsection{Weak limits of weak solutions}
In order to motivate our concept of measure--valued solution,
we consider a family of initial data
satisfying
\begin{equation} \label{MV6}
\begin{split}
\vr_{0, \ep} > 0, \ \intO{ \vr_{0, \ep} } \geq M_0 > 0,\ \ \vt_{0,\ep} &> 0, \ \log \left( \frac{ \vt_{0,\ep}^{c_v} }{\vr_{0, \ep} } \right)  \geq s_0 > - \infty, \\
\intO{ \left[ \frac{1}{2} \vr_{0, \ep}  |\vu_{0, \ep} |^2 + c_v \vr_{0, \ep} \vt_{0, \ep} \right] } &\leq e_0,
\end{split}
\end{equation}
uniformly for $\ep \to 0$.
Suppose that $\{ \vre, \vte, \vue \}_{\ep > 0}$ are the corresponding weak solutions to the Euler system specified through (\ref{MV2}--\ref{MV5}). Our goal is to
identify the cluster point of $\{ \vre, \vte, \vue \}_{\ep > 0}$ for $\ep \to 0$. To this end, we first derive the available {\it a priori} bounds.

\subsubsection{{ A priori} bounds}

To begin, consider $Z \in BC(R)$,
\[
Z' \geq 0, \ Z(s)  \left\{ \begin{array}{lr} < 0 \ \mbox{for}\ s < s_0, \\ \\ =0 \ \mbox{for}\ s \geq s_0. \end{array} \right.
\]
Take $\varphi = \psi(t)$, $\psi \geq 0$ as a test function in the entropy inequality (\ref{MV5}). Using (\ref{MV6}) we deduce, after a straightforward manipulation, that
\[
\intO{  \vre(\tau, \cdot)  Z\left( \log \left( \frac{ \vte^{c_v}(\tau, \cdot) }{\vre (\tau, \cdot) } \right) \right) } \geq 0
\ \mbox{for a.a.} \ \tau \in (0,T);
\]
whence
\begin{equation} \label{MV7}
\log \left( \frac{ \vte^{c_v} }{\vre } \right) \geq s_0
\ \mbox{whenever}\ \vre > 0 \ \Leftrightarrow \ \vre \leq \exp (-s_0) \vte^{c_v} \ \mbox{for a.a. } (\tau, x) \in (0,T) \times \Omega.
\end{equation}

Similarly, we deduce from (\ref{MV4}) and \eqref{MV6} that
\begin{equation} \label{MV8}
\intO{ \left[ \frac{1}{2} \vr_{\ep}  |\vu_{\ep} |^2 + c_v \vr_{\ep} \vt_{\ep} \right](\tau, \cdot) } \leq e_0
\ \mbox{for a.a.}\ \tau \in (0,T).
\end{equation}

Next, in view of (\ref{MV7}),
\[
\vre^{1+\frac1{c_v}}  \leq c( s_0)  c_v \vre \vte;
\]
whence, in accordance with (\ref{MV8}),
\begin{equation} \label{MV9}
\intO{ \vre^{1+\frac1{c_v}} (\tau, \cdot) } \leq c( s_0, e_0) \ \mbox{for a.a.}\ \tau \in (0,T).
\end{equation}

By the same token,
\[
\vre |\log (\vte) |^q \leq \left\{ \begin{array}{lr} c( s_0) \vte^{c_v}| \log(\vte)|^q \leq c(q,s_0)  &\mbox{if}\ \vte \leq 1, \\ \\
\vre \vte  &\mbox{if}\ \vte \geq 1, \end{array} \right.
\]
therefore
\begin{equation} \label{MV10}
\intO{ \vre |\log (\vte)|^q (\tau, \cdot) } \leq c(q, s_0, e_0), \ q \geq 1,  \ \mbox{for a.a.}\ \tau \in (0,T).
\end{equation}

Finally, writing $\vre \vue = \sqrt{\vre} \sqrt{\vre} \vue$, we deduce from (\ref{MV8}) and (\ref{MV9}) that
\begin{equation} \label{MV10extra}
\intO{ |\vre \vue|^p (\tau, \cdot) } \leq c(s_0,e_0) \ \mbox{for a.a.} \ \tau \in (0,T)
\ \mbox{and some}\ p > 1.
\end{equation}

\subsubsection{Young measure}

Unfortunately, the {\it a priori} bounds available are not strong enough to
perform the pointwise limit in the nonlinearities in the weak formulation. Instead we use the characterization of limits of oscillatory
sequences of functions via \emph{Young measures}.
Note that there is an additional  problem as all bounds obtained in the previous section depend on $\vre$. In other words, we have no control over the behavior of
$\vue$, $\vte$ on the (hypothetical) vacuum zone. Consequently, it is more convenient to work with a new set of
state variables - the density $\vr$, the momentum $\vc{m} = \vr \vu$, and the internal energy density $E = c_v \vr \vt$ - the norm of which
is controlled at least in the Lebesgue space $L^1$.

Let
\[
\mathcal{F} = \left\{ [\vr, E, \vc{m} ] \ \Big| \ \vr \in [0, \infty), \ E \in [0, \infty), \ \vc{m} \in R^3 \right\},
\]
denote the new state space.
By virtue of the fundamental theorem on Young measures, see e.g. Ball \cite{BALL2}, there exists a subsequence (not relabeled)
of $\{ \vre, E_\ep \equiv c_v \vre \vte, \vc{m}_\ep \equiv \vre \vue \}_{\ep > 0}$ and a parameterized family of probability measures $\{ Y_{t,x} \}_{(t,x) \in (0,T) \times \Omega}$,
\[
[ (t,x) \mapsto Y_{t,x}] \in L^\infty_{\rm weak-(*)} ((0,T) \times \Omega;
\mathcal{P}(\mathcal{F}))
\]
such that
\begin{equation} \label{MV11}
\left< Y_{t,x}; G(\vr, E, \vc{m}) \right> = \Ov{ G(\vr, E, \vc{m}) }(t,x) \ \mbox{for any}\ G \in C_c(\mathcal{F}) \ \mbox{and a.a.}\ (t,x) \in (0,T) \times \Omega,
\end{equation}
whenever
\begin{equation} \label{MV11a}
G(\vre, c_v \vre \vte, \vre \vue) \to \Ov{ G(\vr, E, \vc{m}) } \ \mbox{weakly-(*) in}\ L^\infty((0,T) \times \Omega).
\end{equation}
The parameterized family of measures $\{ Y_{t,x} \}_{t,x \in (0,T) \times \Omega}$ is called \emph{Young measure} associated to the sequence
$\{ \vre, c_v \vre \vte, \vre \vue \}_{\ep > 0}$.
As a consequence of (\ref{MV7}), we get
\begin{equation} \label{MV12}
{\rm supp} [Y_{t,x}] \subset \left\{ [\vr, E, \vc{m} ] \in \mathcal{F} \ \Big| \ \vr^{1 + c_v} \leq c_v^{-c_v} \exp(-s_0) E^{c_v} \right\}.
\end{equation}

As the nonlinearities appearing in the weak formulation do not in general belong to the class $C_c(\mathcal{F})$ (in the new set of variables $[\vr, E, \vc{m}]$),  validity of
(\ref{MV11a}) must be extended to a larger class of functions.
If $G \in C (\mathcal{F})$ is such that
\begin{equation} \label{MV13}
\int_0^T \intO{ |G(\vre, c_v \vre \vte, \vre \vue)| } \leq c \ \mbox{uniformly for}\ \ep \to 0,
\end{equation}
then $G$ in $Y_{t,x}$ integrable for a.a. $(t, x) \in (0,T) \times \Omega$ and
\[
\left[ (t,x) \mapsto \left< Y_{t,x}; G(\vr, E, \vc{m}) \right> \right] \in L^1((0,T) \times \Omega).
\]
The function $\left[ (t,x) \mapsto \left< Y_{t,x}; G(\vr, E, \vc{m} ) \right> \right]$ can be identified with the so-called \emph{biting limit} of
the family $\{ G(\vre, c_v \vre \vte, \vre \vue) \}_{\ep > 0}$, see Ball and Murat \cite{BAMU}. Finally, the same holds for
any $G : \mathcal{F} \to R \cup \{ \infty \}$ satisfying (\ref{MV13}) and such that there exists a sequence $G_m \in C_c(\mathcal{F})$, $G_m \nearrow G$ in $\mathcal{F}$.

If (\ref{MV13}) holds, we have
\[
G(\vre, c_v \vre \vte, \vre \vue) \to \Ov{G(\vr, E, \vc{m})} \ \mbox{weakly-(*) in} \ \mathcal{M} ([0,T] \times \Omega),
\]
for a suitable subsequence. Here, the singular part of the limit measure reflects possible concentrations in  $\{ \vre, \vre \vte, \vre \vue \}_{\ep > 0}$.

\begin{Remark} \label{RMV1}

Note carefully that the Young measure $\left[ (t,x) \mapsto \left< Y_{t,x}; G(\vr, E, \vc{m} ) \right> \right]$ is a parameterized family of
non-negative measures acting on the phase space $\mathcal{F}$ while $\Ov{G(\vr, E, \vc{m})}$ is a signed measure on the physical space $[0,T] \times \Omega$.

\end{Remark}

The difference
\[
\mu_G \equiv \Ov{G(\vr, E, \vc{m})} - \left[ (t,x) \mapsto \left< Y_{t,x}; G(\vr, E, \vc{m}) \right> \right] \in \mathcal{M} ([0,T] \times \Omega),
\]
is called \emph{concentration defect measure}. It vanishes whenever the family
$\{ G(\vre, c_v \vre \vte, \vre \vue) \}_{\ep > 0}$ is equi-integrable (weakly precompact) in $L^1((0,T) \times \Omega)$.

We claim the following result proved in \cite[Lemma 2.1]{FGSWW1}.

\begin{Lemma} \label{LMV1}

Let
\[
|G(\vr, E, \vc{m} )| \leq F(\vr, E, \vc{m}) \ \mbox{for all} \ (\vr, E, \vc{m}) \in \mathcal{F}.
\]

Then
\[
\left| \Ov{G(\vr, E, \vc{m})} - \left< Y_{t,x}; G(\vr, E, \vc{m}) \right> \right|
\leq \Ov{F(\vr, E, \vc{m})} - \left< Y_{t,x}; F(\vr, E, \vc{m}) \right> \equiv \mu_F
\ \mbox{in}\ \mathcal{M} ([0,T] \times \Omega).
\]

\end{Lemma}

\subsubsection{The limit $\ep \to 0$}

We are ready to perform the limit for $\ep \to 0$ in the weak formulation (\ref{MV2}--\ref{MV5}).

\medskip

{\bf Step 1}

In view of the uniform bounds (\ref{MV9}) and  (\ref{MV10extra}), we can let $\ep \to 0$ in (\ref{MV2}) obtaining
\begin{equation} \label{MV14}
\int_0^T \intO{ \left[ \left< Y_{t,x}; \vr \right> \partial_t \varphi + \left< Y_{t,x}; \vc{m} \right> \cdot \Grad \varphi \right] } \dt = - \intO{
\left< Y_{0,x}; \vr \right> \varphi(0, \cdot) },
\end{equation}
for any $\varphi \in \DC([0,T) \times \Omega)$, where $Y_{0,x}$ is the Young measure generated by the initial data. Note that, in accordance with
(\ref{MV9}) and (\ref{MV10extra}), the
families $\{ \vre \}_{\ep > 0}$, $\{ \vre \vue \}_{\ep > 0}$ are equi-integrable and therefore concentrations do not occur. Finally, we deduce from
(\ref{MV14}) that
\begin{equation} \label{MV15}
\left[ \intO{
\left< Y_{t,x}; \vr \right> \varphi } \right]_{t = 0}^{t = \tau} =
\int_0^\tau \intO{ \left[ \left< Y_{t,x}; \vr \right> \partial_t \varphi + \left< Y_{t,x}; \vc{m} \right> \cdot \Grad \varphi \right] } \dt,
\end{equation}
for a.a. $\tau \in (0,T)$ and for any $\varphi \in C^1 ([0,T] \times \Omega)$.

\begin{Remark} \label{RMV2}

Relation (\ref{MV15}) can be justified for \emph{any} $\tau \in [0,T)$, however, this is not needed for future analysis.

\end{Remark}

\medskip

\medskip

{\bf Step 2}

Keeping in mind the uniform energy bound (\ref{MV8}),
we consider $\varphi = \psi (t)$ as a test function in the total energy balance (\ref{MV4}) obtaining
\[
\begin{split}
&\int_0^T \intO{ \left< Y_{t,x} ;  \frac{1}{2} \frac{ |\vc{m} |^2 }{\vr} + E  \right> \partial_t \psi
 } \dt \\&+ \int_0^T \left< \Ov{ \left(\frac{1}{2} \frac{ |\vc{m} |^2 }{\vr} + E\right)}(t,\cdot) ; \Omega \right> \partial_t \psi \ \dt  -
 \int_0^T \intO{ \left< Y_{t,x} ;  \frac{1}{2} \frac{ |\vc{m} |^2 }{\vr} + E  \right> \partial_t \psi } \dt
\\ &= - \psi(0) \intO{ \left< Y_{0,x} ; \left( \frac{1}{2} \frac{ |\vc{m} |^2 }{\vr} + E \right) \right>  }.
\end{split}
\]

\begin{Remark} \label{RMvW1}

Note that the function
\[
[ \vr , \vc{m} ] \in {\rm int}[\mathcal{F}] \mapsto \frac{|\vc{m}|^2}{\vr},
\]
extended to be $0$ whenever $\vc{m} = 0$ and $\infty$ if $\vr = 0$, $\vc{m} \ne 0$ is a convex lower semi-continuous in $\mathcal{F}$;
whence $Y_{t,x}$ measurable.

\end{Remark}

Thus we may infer that
\begin{equation} \label{MV19}
\left[  \intO{ \left< Y_{t,x} ;  \frac{1}{2} \frac{ |\vc{m} |^2 }{\vr} + E  \right> } \right]_{t=0}^{t = \tau} +
\mathcal{D}(\tau) = 0,
\end{equation}
for a.a. $\tau \in (0,T)$,
with non-negative $\mathcal{D} \in L^\infty(0,T)$, 
\begin{equation} \label{MV20}
\mathcal{D} (\tau) =
\left< \Ov{\left( \frac{1}{2} \frac{ |\vc{m} |^2 }{\vr} + E  \right) }(\tau, \cdot)  ; \Omega  \right>    -  \intO{ \left< Y_{\tau,x} ; \frac{1}{2} \frac{ |\vc{m} |^2 }{\vr} + E  \right>  },  
\end{equation}
for a.a. $\tau \in (0,T)$. The quantity $\mathcal{D}$ will be termed \emph{dissipation defect}.

\begin{Remark} \label{RVacbis}

Recall that, in view of (\ref{MV8}), we have
\[
\Ov{\left( \frac{1}{2} \frac{ |\vc{m} |^2 }{\vr} + E  \right) } \in L^\infty(0,T; \mathcal{M}(\Omega)).
\]

\end{Remark}

\begin{Remark} \label{RVac}

It follows from (\ref{MV19}) that
\begin{equation} \label{MV12a}
{\rm supp}[Y_{t,x}] \cap \{ \vr = 0 , \ \vc{m} \ne 0 \} = \emptyset.
\end{equation}
In particular, if the measure $Y_{t,x}$ charges vacuum zone it must be only within the hyperplane $\vc{m}  = 0$.

\end{Remark}

\medskip

{\bf Step 3}

Similarly, we deduce from the momentum equation (\ref{MV3}) that
\begin{equation} \label{MV16}
\begin{split}
&\left[ \intO{ \left< Y_{t,x}; \vc{m} \right> \cdot \vcg{\varphi} }
\right]_{t = 0}^{t = \tau} \\ &= \int_0^\tau \intO{ \left[ \left< Y_{t,x}; \vc{m} \right> \cdot \partial_t\vcg{\varphi} +
\left< Y_{t,x}; \frac{\vc{m} \otimes \vc{m} }{\vr} \right> : \Grad \vcg{\varphi} + \left< Y_{t,x} ; p \right> \Div \vcg{\varphi} \right] } \dt\\
&+ \int_0^\tau \Grad \vcg{\varphi }:{\rm d}\mu_{\mathcal{R}},
\end{split}
\end{equation}
for a.a. $\tau \in (0,T)$ and for any $\vcg{\varphi} \in C^1([0,T] \times \Omega; R^3)$,
with the concentration error
\[
\mu_{\mathcal{R}} = \Ov{ \left( \frac{\vc{m} \otimes \vc{m} }{\vr} \right)} - \left< Y_{t,x}; \frac{\vc{m} \otimes \vc{m} }{\vr} \right> +
\Ov{ p } \mathbb{I} - \left< Y_{t,x}; p \right> \mathbb{I} \in L^\infty(0,T; \mathcal{M} (\Omega; R^{3 \times 3} )).
\]
Note that
$p = \frac{1}{c_v} E$, and that, in view of Remark \ref{RVac}, it is enough to set $\frac{\vc{m} \otimes \vc{m}}{\vr} = 0$ whenever $\vc{m} = 0$.

In accordance with Lemma \ref{LMV1}, we have an important relation
between the concentration error $\mu_{\mathcal{R}}$ in (\ref{MV16}) and the dissipation defect $\mathcal{D}$, namely
\begin{equation} \label{MV21}
\left| \int_0^\tau \Grad \vcg{\varphi} : {\rm d} \mu_{\mathcal{R}} \right| \leq
\| \Grad \vcg{\varphi} \|_{C([0,\tau] \times \Omega)} \int_0^\tau \mathcal{D}(t) \ \dt
\ \mbox{for a.a.} \ \tau \in (0,T).
\end{equation}
\medskip

{\bf Step 4}

Finally, the entropy balance (\ref{MV5}) gives rise to
\begin{equation} \label{MV17}
\begin{split}
&\left[
\intO{ \left< Y_{t,x}; \vr Z \left( s \right) \right> \varphi  } \right]_{t = 0}^{t = \tau}\\
&\geq \int_0^\tau \intO{ \left[ \left< Y_{t,x} ; \vr Z\left(s  \right) \right>  \partial_t \varphi +
\left< Y_{t,x} ; \vr Z \left( s\right) \vc{m} \right> \cdot \Grad \varphi  \right] } \dt
\end{split}
\end{equation}
for a.a. $\tau \in (0,T)$, any $\varphi \in C^1([0,T] \times \Omega)$, $\varphi \geq 0$, $Z \in BC(R)$, $Z' \geq 0$,
where
\[
s =s(\vr, E) = \log \left( \frac{E^{c_v}}{\vr^{c_v + 1} } \right).
\]

The couple $\{ Y_{t,x}, \mathcal{D} \}$ satisfying (\ref{MV15}), (\ref{MV19}), (\ref{MV16}), (\ref{MV21}), and (\ref{MV17})
represents a \emph{dissipative measure-valued solution} of the complete Euler system (\ref{i1}--\ref{i4}).

\subsection{Dissipative measure-valued solutions}

Motivated by the previous discussion, we introduce the concept of a dissipative measure-valued solution to the Euler system (\ref{i1}--\ref{i4})
for general constitutive relations. Although motivated by the preceding section, the measure-valued solutions
introduced below represent an object formally independent of any approximation procedure, in particular they may not be a limit of a family
of weak solutions.

In addition to Gibbs' equation (\ref{i5}), we assume the \emph{hypothesis of thermodynamic stability},
\begin{equation} \label{WS1}
\frac{\partial p(\vr, \vt)}{\partial \vr} > 0,\ \frac{\partial e(\vr, \vt)}{\partial \vt} > 0 \ \mbox{for all}\ \vr, \vt > 0.
\end{equation}
In particular, any function $G = G(\vr, \vt, \vu)$ can be identified with a function of variables $[\vr, E = \vr e(\vr,\vt), \vc{m} = \vr \vu]$
as
\[
G(\vr, \vt, \vu) = G \left(\vr, \vt(\vr, E), \frac{\vc{m} }{\vr} \right) \ \mbox{for all}\ \vr > 0, \ \vt > 0, \ \vu \in R^3.
\]
We simply write $G(\vr, E, \vc{m})$ as the case may be.

\begin{Remark} \label{bfrr2}
The former condition in (\ref{WS1}) means that the compressibility of the gas is positive while the
latter expresses positivity of the specific heat at constant volume. They may can be rephrased as
convexity of the internal energy $e$ as a function of the entropy $s$ and the specific volume $\frac{1}{\vr}$,
see Bechtel, Rooney and Forest \cite{BEROFO}.
\end{Remark}

\begin{Definition} \label{DMV1} {\bf [Dissipative measure-valued solution]}

{\it
A family of probability measures $\{ Y_{t,x} \}_{(t,x) \in (0,T) \times \Omega}$,
\[
(t,x) \mapsto Y_{t,x} \in L^\infty_{\rm weak-(*)} ((0,T) \times \Omega; \mathcal{P}(\mathcal{F})),\
\]
and the dissipation defect $\mathcal{D} \in L^\infty(0,T)$ represent a \emph{dissipative measure-valued solution} of the Euler system
(\ref{i1}--\ref{i4})
with the initial data $Y_{0,x}$ if:

$\bullet$
\begin{equation} \label{MV22}
\begin{split}
\left[ \intO{
\left< Y_{t,x}; \vr \right> \varphi } \right]_{t = 0}^{t = \tau} =
\int_0^\tau \intO{ \left[ \left< Y_{t,x}; \vr \right> \partial_t \varphi + \left< Y_{t,x}; \vc{m} \right> \cdot \Grad \varphi \right] } \dt,
\end{split}
\end{equation}
for a.a. $\tau \in (0,T)$ and for any $\varphi \in C^1 ([0,T] \times \Omega)$;

$\bullet$
\begin{equation} \label{MV23}
\begin{split}
&\left[ \intO{ \left< Y_{t,x}; \vc{m} \right> \cdot \vcg{\varphi} }
\right]_{t = 0}^{t = \tau} \\ &= \int_0^\tau \intO{ \left[ \left< Y_{t,x}; \vc{m} \right> \cdot \partial_t\vcg{\varphi} +
\left< Y_{t,x}; \frac{\vc{m} \otimes \vc{m} }{\vr}  \right> : \Grad \vcg{\varphi} + \left< Y_{t,x} ; p(\vr,E) \right> \Div \vcg{\varphi} \right] } \dt\\
&+ \int_0^\tau \Grad \vcg{\varphi }: {\rm d}\mu_{\mathcal{R}},
\end{split}
\end{equation}
for a.a. $\tau \in (0,T)$ and for any $\vcg{\varphi} \in C^1([0,T] \times \Omega; R^3)$;

$\bullet$
\begin{equation} \label{MV24}
\begin{split}
\bigg[
\int_{\Omega}  \big< Y_{t,x}; \vr Z  &    \left( s(\rho, E) \right) \big> \varphi  \dx \bigg]_{t = 0}^{t = \tau} 
\\ & \geq \int_0^\tau \intO{ \Big[ \left< Y_{t,x} ; \vr Z\left( s(\rho, E) \right) \right>  \partial_t \varphi +
\left< Y_{t,x} ; Z \left( s(\rho, E) \right) \vc{m} \right> \cdot \Grad \varphi  \Big] } \dt,
\end{split}
\end{equation}
for a.a. $\tau \in (0,T)$, any $\varphi \in C^1([0,T] \times \Omega)$, $\varphi \geq 0$, $Z \in BC(R)$, $Z' \geq 0$;

$\bullet$
\begin{equation} \label{MV25}
\left[  \intO{ \left< Y_{t,x} ;  \frac{1}{2} \frac{|\vc{m}|^2 }{\vr} + E  \right> } \right]_{t=0}^{t = \tau} +
\mathcal{D}(\tau) = 0,
\end{equation}
where the dissipation defect $\mathcal{D}$ dominates the signed measure
\[
\mu_{\mathcal{R}} \in \mathcal{M}([0,T] \times \Omega; R^{3 \times 3}),
\]
specifically,
\begin{equation} \label{MV27}
 \| \mu_{\mathcal{R}}
\|_{\mathcal{M}([0,\tau) \times \Omega; R^{3 \times 3})} \leq c \int_0^\tau \mathcal{D}(t) \ \dt,
\end{equation}
for a.a. $\tau \in (0,T)$. }
\end{Definition}

As already pointed out, the dissipative measure-valued solutions are designed to retain the \emph{minimal} piece of information inherited from the original Euler system
in the course of some limit process. They may be seen as limits of families of weak solutions or their numerical approximations, cf. \cite{FGSWW1}.
Notably, as shown in the next section, the dissipative measure-valued solutions comply
with the {\it weak--strong uniqueness principle}. In this context, relation (\ref{MV27}) plays the crucial role.

\section{Weak--strong uniqueness}
\label{WS}

Our ultimate goal is to show the main result of the present paper, namely, a dissipative measure valued solution and a strong solution starting from the same initial data coincide as long as the latter exists. In addition to the natural physical principles encoded in (\ref{i5}), (\ref{WS1}), we shall need a purely technical hypothesis
\begin{equation} \label{WS1a}
|p(\vr, \vt)| \aleq (1 + \vr + \vr e(\vr, \vt) + \vr |s(\vr, \vt)| ).
\end{equation}

Note that (\ref{WS1a}) is satisfied for a large family of gases for which $p \approx \vr \vt$ including
the perfect gas studied in Section \ref{MV}.
Here and hereafter, the symbol $a \aleq b$ means $a \leq c b$ for a certain constant $c > 0$.

\subsection{Relative energy}

Let
\begin{equation} \label{WS6}
r \in C^1([0,T] \times \Omega), \ r > 0, \ \Theta \in C^1([0,T] \times \Omega), \ \Theta > 0, \ \vc{U} \in C^1([0,T] \times \Omega; R^3),
\end{equation}
be given. Following \cite{FeiNov10}, we introduce the \emph{ballistic free energy}
\[
H_\Theta(\vr, \vt) = \vr e(\vr, \vt) - \Theta \vr s(\vr, \vt),
\]
and the \emph{relative energy}
\[
\mathcal{E}_Z \left(\vr, \vt, \vu \ \Big| r, \Theta, \vc{U} \right) =
\frac{1}{2} \vr |\vu - \vc{U}|^2 + \vr e(\vr, \vt) - \Theta \vr Z(s(\vr, \vt)) -
\frac{ \partial H_\Theta(r, \Theta) }{\partial \vr} (\vr - r) - H_\Theta (r, \Theta).
\]
The relative energy can be written in the new variables $[\vr, E, \vc{m}]$ as
\begin{equation} \label{WS2}
\mathcal{E}_Z \left(\vr, E, \vc{m} \ \Big| r, \Theta, \vc{U} \right) =
\frac{1}{2} \vr \left| \frac{\vc{m}}{\vr} - \vc{U} \right|^2 + E - \Theta \vr Z(s(\vr, E)) -
\frac{ \partial H_\Theta(r, \Theta) }{\partial \vr} (\vr - r) - H_\Theta (r, \Theta).
\end{equation}
In contrast with \cite{FeiNov10}, the relative entropy functional depends also on the cut-off function appearing in the entropy inequality
(\ref{MV24}). The specific shape of $Z$ will be fixed below.

\begin{Remark} \label{bfrr4}
Notation in (\ref{WS2}) is slightly inconsistent as $\frac{\partial H_\Theta(r, \Theta)}{\partial \vr}$ still denotes the derivative
with respect to $\vr$ of the function $H_\Theta(\vr, \vt)$ considered in the ``old'' variables $(\vr, \vt)$ rather than $(\vr, E)$.
We still believe this is convenient as the ``test functions'' $r$ and $\Theta$ are designed to mimick the density and the absolute temperature
of the strong solution.
\end{Remark}

\subsubsection{Relative energy inequality}

Using the abstract formulation (\ref{MV22}--\ref{MV25}) we derive a functional relation
\begin{equation} \label{WS3}
\begin{split}
&\left[ \intO{ \left< Y_{t,x} ; \mathcal{E}_Z \left(\vr, E, \vc{m} \ \Big| r, \Theta, \vc{U} \right) \right> } \right]_{t = 0}^{t = \tau} =
\left[ \intO{ \left< Y_{t,x} ; \frac{1}{2} \frac{|\vc{m}|^2}{\vr} + E \right> } \right]_{t = 0}^{t = \tau}\\
& - \left[ \intO{ \left< Y_{t,x}; \vc{m} \right> \cdot \vc{U} } \right]_{t = 0}^{t = \tau} +
\left[ \intO{ \left< Y_{t,x}; \vr \right> \left( \frac{1}{2} |\vc{U}|^2 - \frac{ \partial H_\Theta(r, \Theta) }{\partial \vr} \right)  } \right]_{t = 0}^{t = \tau}\\
& - \left[ \intO{ \left< Y_{t,x} ; \vr Z(s(\vr, E)) \right> \Theta } \right]_{t = 0}^{t = \tau}
+ \left[ \intO{ \frac{ \partial H_\Theta(r, \Theta) }{\partial \vr} r - H_\Theta(r, \Theta) } \right]_{t = 0}^{t = \tau}\\
&= - \mathcal{D}(\tau)  - \left[ \intO{ \left< Y_{t,x}; \vc{m} \right> \cdot \vc{U} } \right]_{t = 0}^{t = \tau} +
\left[ \intO{ \left< Y_{t,x}; \vr \right> \left( \frac{1}{2} |\vc{U}|^2 - \frac{ \partial H_\Theta(r, \Theta) }{\partial \vr} \right)  } \right]_{t = 0}^{t = \tau}\\
& - \left[ \intO{ \left< Y_{t,x} ; \vr Z(s(\vr, E)) \right> \Theta } \right]_{t = 0}^{t = \tau}
+ \left[ \intO{ \frac{ \partial H_\Theta(r, \Theta) }{\partial \vr} r - H_\Theta(r, \Theta) } \right]_{t = 0}^{t = \tau}.
\end{split}
\end{equation}

Furthermore, using the entropy inequality (\ref{MV24}), we get
\begin{equation} \label{WS4}
\begin{split}
&\left[ \intO{ \left< Y_{t,x} ; \mathcal{E}_Z \left(\vr, E, \vc{m} \ \Big| r, \Theta, \vc{U} \right) \right> } \right]_{t = 0}^{t = \tau}
+ \mathcal{D}(\tau)
\\
&\leq
- \int_0^\tau \intO{ \left[ \left< Y_{t,x} ; \vr Z\left( s(\vr,E) \right) \right>  \partial_t \Theta +
\left< Y_{t,x} ; Z \left( s(\vr, E) \right) \vc{m} \right> \cdot \Grad \Theta  \right] } \dt\\
& - \left[ \intO{ \left< Y_{t,x}, \vc{m} \right> \cdot \vc{U} } \right]_{t = 0}^{t = \tau} +
\left[ \intO{ \left< Y_{t,x}; \vr \right> \left( \frac{1}{2} |\vc{U}|^2 - \frac{ \partial H_\Theta(r, \Theta) }{\partial \vr} \right)  } \right]_{t = 0}^{t = \tau}\\
&
+ \left[ \intO{ \frac{ \partial H_\Theta(r, \Theta) }{\partial \vr} r - H_\Theta(r, \Theta) } \right]_{t = 0}^{t = \tau}.
\end{split}
\end{equation}

The advantage of (\ref{WS4}) is that all integrals on its right-hand side can be expressed by means of (\ref{MV22}) and (\ref{MV23}). Thus, repeating the arguments of \cite[Section 3]{FeiNov10}, we obtain
\begin{equation} \label{WS5}
\begin{split}
&\left[ \intO{ \left< Y_{t,x} ; \mathcal{E}_Z \left(\vr, E, \vc{m} \ \Big| r, \Theta, \vc{U} \right) \right> } \right]_{t = 0}^{t = \tau}
+ \mathcal{D}(\tau)
\\
&\leq
- \int_0^\tau \intO{ \left[ \left< Y_{t,x} ; \vr Z\left( s(\vr,E) \right) \right>  \partial_t \Theta +
\left< Y_{t,x} ;  Z \left( s(\vr,E)  \right) \vc{m} \right> \cdot \Grad \Theta  \right] } \dt\\
&+ \int_0^\tau \intO{ \left[ \left< Y_{t,x}; \vr \right> s(r, \Theta) \partial_t \Theta + \left< Y_{t,x}; \vc{m} \right> \cdot s(r, \Theta) \Grad \Theta    \right] } \dt
\\
&+ \int_0^\tau \intO{ \left[ \left< Y_{t,x}; \vr \vc{U} - \vc{m} \right> \cdot \partial_t \vc{U} +
\left< Y_{t,x}; (\vr \vc{U} - \vc{m} ) \otimes \frac{\vc{m}}\vr \right> : \Grad \vc{U} - \left< Y_{t,x}; p(\vr, E) \right> \Div \vc{U}  \right] } \dt \\
&+ \int_0^\tau \intO{ \left[ \left< Y_{t,x} ; r - \vr \right> \frac{1}{r} \partial_t p(r, \Theta) -
\left< Y_{t,x} ; \vc{m} \right> \cdot \frac{1}{r} \Grad p(r, \Theta)    \right] } \ \dt\\
&- \int_0^\tau \Grad \vc{U}:{\rm d}\mu_{\mathcal{R}}.
\end{split}
\end{equation}

The relation (\ref{WS5}) holds for any dissipative measure-valued solution of the Euler system and any trio of smooth test functions satisfying
(\ref{WS6}). It can be seen as a measure-valued variant of the relative energy inequality derived in \cite{FeiNov10}.

\begin{Remark} \label{bfrrr14}
The fact that (\ref{WS5}) holds for \emph{any} trio of ``test functions''
$[r, \Theta, \vc{U}]$ is important in view of possible future applications, cf. Section \ref{C}.
\end{Remark}

\subsection{Weak-strong uniqueness in the class of measure-valued solutions}

Suppose that $[r, \Theta, \vc{U}]$ is a strong solution of the Euler system (\ref{i1}--\ref{i3}) starting from the initial data
$[r_0, \Theta_0, \vc{U}_0]$ belonging to the class (\ref{WS6}). We fix a compact set $K \subset (0, \infty)^2$ containing the trajectories
$\cup_{t \in [0,T], x \in \Omega} [r(t,x), \Theta (t,x)]$ and its image $\tilde K \subset (0, \infty)^2$ in the new phase variables
\[
(\vr, \vt) \mapsto [\vr, \vr e(\vr, \vt)] : (0,\infty)^2 \to (0, \infty)^2.
\]
Finally, we consider a function $\Phi (\vr, E)$,
\[
\Phi \in \DC (0,\infty)^2, \ 0 \leq \Phi \leq 1,\  \Phi|_{\mathcal{U}} = 1, \ \mbox{where}\ \mathcal{U} \ \mbox{is an open neighborhood of} \ \tilde K
\ \mbox{in}\ (0,\infty)^2.
\]
For a measurable function $G(\vr, E, \vc{m})$, we set
\[
G = G_{\rm ess} + G_{{\rm res}} , \ G_{\rm ess} = \Phi (\vr, E) G(\vr, E, \vc{m}),\
G_{\rm res} = (1 - \Phi (\vr, E)) G(\vr, E, \vc{m}).
\]
The idea, borrowed from \cite{FENO6}, is that the ``essential part'' $G_{\rm ess}$ describes the behavior of the non-linearity in the non-degenerate 
area where both $\vr$ and $\vt$ are bounded below and above, while    
the ``residual part'' $G_{\rm res}$ captures the behavior in the singular regime $\vr, \vt \to 0$ or/and $\vr, \vt \to \infty$. 

Finally, we consider $Z = Z_{a,b} \in BC(R)$, $-\infty \leq a < b \leq \infty$,
\[
Z_{a,b} (s) = \left\{ \begin{array}{l} a \ \mbox{for} \ s < a, \\ s \ \mbox{for}\ s \in [a,b], \\
b \ \mbox{for}\ s \geq b, \end{array} \right. 
\]
and fix $a,b$ finite in such a way that
\begin{equation} \label{WS5b}
[ Z_{a,b} (s(\vr, E) ) ]_{\rm ess} = \Phi (\vr, E) Z_{a,b} (s(\vr, E)) = \Phi (\vr, E)s(\vr, E) = [s (\vr, E)]_{\rm ess}.
\end{equation}

\subsubsection{Initial data}

We consider a dissipative measure valued solution $\{ Y_{t,x}, \mathcal{D} \}$ such that its initial value coincides with
$[r_0, \Theta_0, \vc{U}_0]$, meaning
\[
Y_{0,x} = \delta_{[ r_0(x), r_0 e(r_0, \Theta_0)(x), r_0 \vc{U}_0 (x) ]} \ \mbox{for a.a.}\ x \in \Omega,
\]
where $\delta_Y$ denotes the Dirac distribution supported at $Y$.
Accordingly,
\[
\intO{ \left< Y_{0,x} ; \mathcal{E}_Z \left( \vr, E, \vc{m} \ \Big| r_0(x), \Theta_0(x), \vc{U}_0(x) \right) \right> } = 0.
\]

Taking $\varphi = 1$ in (\ref{MV24}) we get
\[
\intO{ \left< Y_{\tau,x}; \vr Z (s(\vr, E)) \right> } \geq \intO{ r_0 Z (s(r_0, \Theta_0 )) } \ \mbox{for a.a.} \ \tau \in (0,T).
\]
As the initial data are regular, we deduce that there exists $a \in R$ such that
\[
\intO{ \left< Y_{\tau,x}; \vr Z (s(\vr,E)) \right> } = 0 \ \mbox{whenever} \ Z \leq 0, \ Z(s) = 0\ \mbox{for all} \ s \geq a.
\]
Consequently, we obtain
\[
\intO{ \left< Y_{\tau,x}; \vr Z_{a,b} (s(\vr,E)) \right> } = \intO{ \left< Y_{\tau,x}; \vr Z_{-\infty,b} (s(\vr,E)) \right> },
\]
in particular
\begin{equation} \label{WS5a}
- \intO{ \left< Y_{\tau,x}; \vr Z_{a,b} (s(\vr,E)) \right> } = - \intO{ \left< Y_{\tau,x}; \vr Z_{-\infty,b} (s(\vr,E)) \right> }
\geq - \intO{ \left< Y_{\tau,x}; \vr s(\vr,E) \right> }.
\end{equation}

Thus introducing a new relative energy
\[
\mathcal{E} \left(\vr, E, \vc{m} \ \Big| r, \Theta, \vc{U} \right) =
\frac{1}{2} \vr \left| \frac{\vc{m}}{\vr} - \vc{U} \right|^2 + E - \Theta \vr s(\rho,E)  -
\frac{ \partial H_\Theta(r, \Theta) }{\partial \vr} (\vr - r) - H_\Theta (r, \Theta),
\]
and going back to (\ref{WS5}), we obtain
\begin{equation} \label{WS7}
\begin{split}
& \intO{ \left< Y_{\tau,x} ; \mathcal{E} \left(\vr, E, \vc{m} \ \Big| r, \Theta, \vc{U} \right) \right> }
+ \mathcal{D}(\tau)
\\
&\leq
- \int_0^\tau \intO{ \left[ \left< Y_{t,x} ; \vr Z\left( s(\vr,E)  \right) \right>  \partial_t \Theta +
\left< Y_{t,x} ; Z \left( s(\vr,E) \right) \vc{m}  \right> \cdot \Grad \Theta  \right] } \dt\\
&+ \int_0^\tau \intO{ \left[ \left< Y_{t,x}; \vr \right> s(r, \Theta) \partial_t \Theta + \left< Y_{t,x}; \vc{m} \right> \cdot s(r, \Theta) \Grad \Theta    \right] } \dt
\\
&+ \int_0^\tau \intO{ \left[ \left< Y_{t,x}; \vr \vc{U} - \vc{m}  \right> \cdot \partial_t \vc{U} +
\left< Y_{t,x}; (\vr \vc{U} - \vc{m}) \otimes \frac{\vc{m}}{\vr} \right> : \Grad \vc{U} - \left< Y_{t,x}; p(\vr,E) \right> \Div \vc{U}  \right] } \dt \\
&+ \int_0^\tau \intO{ \left[ \left< Y_{t,x} ; r - \vr \right> \frac{1}{r} \partial_t p(r, \Theta) -
\left< Y_{t,x} ; \vc{m} \right> \cdot \frac{1}{r} \Grad p(r, \Theta)    \right] } \ \dt\\
&- \int_0^\tau \Grad \vc{U}:{\rm d} \mu_{\mathcal{R}} \ \dt  \ \mbox{for a.a.}\ \tau \in (0,T),
\end{split}
\end{equation}
with some fixed $Z = Z_{a,b}$, $a<b$ finite.

\subsubsection{A Gronwall type argument}

Our ultimate goal is to show that the right-hand side of (\ref{WS7}) can be absorbed by the time average of the left-hand side.
Thus, by means of the standard Gronwall argument, the left hand must vanish identically in $(0,T)$.
To this end, we recall the coercivity properties of $\mathcal{E}$ proved in \cite[Chapter 3, Proposition 3.2]{FENO6},
\begin{equation} \label{WS8}
\begin{split}
\mathcal{E} &\left(\vr, E, \vc{m} \ \Big| r, \Theta, \vc{U} \right) \\ &\ageq
\left[ | \vr - r |^2 + |E - r e(r, \Theta) |^2 + \left| \frac{\vc{m}}{\vr} - \vc{U} \right|^2 \right]_{\rm ess} +
\left[ 1 + \vr + \vr |s(\vr,E)| + E +  \frac{|\vc{m}|^2}{\vr} \right]_{\rm res}.
\end{split}
\end{equation}

\medskip

{\bf Step 1}

We first use
(\ref{MV27}) to observe that
\[
\left| \int_0^\tau \Grad \vc{U}:{\rm d} \mu_{\mathcal{R}} \ \dt \right| \aleq \int_0^\tau \mathcal{D}(t) \ \dt.
\]

Next, write 
\[
\begin{split}
\intO{ \left< Y_{t,x}; (\vr \vc{U} - \vc{m} ) \otimes \frac{\vc{m}}{\vr} \right> : \Grad \vc{U} }
&=  \intO{ \left< Y_{t,x}; \vr \vc{U} - \vc{m}  \right> \cdot \vc{U} \cdot \Grad \vc{U} }\\
&+  \intO{ \left< Y_{t,x}; \vr \left(\vc{U} - \frac{\vc{m}}{\vr} \right) \otimes \left( \frac{\vc{m}}{\vr} - \vc{U} \right) \right> : \Grad \vc{U} },
\end{split}
\]
where the right integral is controlled be the left-hand side of (\ref{WS7}).

Consequently, as $[r, \Theta, \vc{U}]$ solve the Euler system
(\ref{i1}--\ref{i3}), inequality (\ref{WS7}) reduces to 
\begin{equation} \label{WS9}
\begin{split}
& \intO{ \left< Y_{\tau,x} ; \mathcal{E} \left(\vr, E, \vc{m} \ \Big| r, \Theta, \vc{U} \right) \right> }
+ \mathcal{D}(\tau)
\\
&\aleq
- \int_0^\tau \intO{ \left[ \left< Y_{t,x} ; \vr Z\left( s(\vr, E)  \right) \right>  \partial_t \Theta +
\left< Y_{t,x} ; Z \left( s(\vr,E)  \right) \vc{m} \right> \cdot \Grad \Theta  \right] } \dt\\
&+ \int_0^\tau \intO{ \left[ \left< Y_{t,x}; \vr \right> s(r, \Theta) \partial_t \Theta + \left< Y_{t,x}; \vc{m} \right> \cdot s(r, \Theta) \Grad \Theta    \right] } \dt
\\
&+ \int_0^\tau \intO{ \left[p(r, \Theta) \Div \vc{U} - \left< Y_{t,x}; p(\vr, E) \right> \Div \vc{U}  \right] } \dt \\
&+ \int_0^\tau \intO{ \left[ \left< Y_{t,x} ; r - \vr \right> \frac{1}{r} \partial_t p(r, \Theta) -
\left< Y_{t,x} ; \vr \vc{U} \right> \cdot \frac{1}{r} \Grad p(r, \Theta) - p(r, \Theta) \Div \vc{U}   \right] } \ \dt\\
& + \int_0^\tau \left[ \intO{ \left< Y_{t,x} ; \mathcal{E} \left(\vr, E, \vc{m} \ \Big| r, \Theta, \vc{U} \right) \right> }
+ \mathcal{D}(t) \right] \dt \
\mbox{for a.a.}\ \tau \in (0,T).
\end{split}
\end{equation}

\medskip

{\bf Step 2}

Keeping in mind (\ref{WS5b}) we may rewrite 
\[
\begin{split}
- &\int_0^\tau \intO{ \left[ \left< Y_{t,x} ; \vr Z\left( s(\vr, E) \right) \right>  \partial_t \Theta +
\left< Y_{t,x} ; Z \left( s(\vr, E) \right) \vc{m} \right> \cdot \Grad \Theta  \right] } \dt\\
+ &\int_0^\tau \intO{ \left[ \left< Y_{t,x}; \vr \right> s(r, \Theta) \partial_t \Theta + \left< Y_{t,x}; \vc{m} \right> \cdot s(r, \Theta) \Grad \Theta    \right] } \dt\\
=& - \int_0^\tau \intO{ \left[ \left< Y_{t,x} ; \left[ \vr Z\left( s(\vr,E) \right) \right]_{\rm ess} \right>  \partial_t \Theta +
\left< Y_{t,x} ; \left[ Z \left( s(\vr, E) \right) \vc{m} \right]_{\rm ess} \right> \cdot \Grad \Theta  \right] } \dt\\
+ &\int_0^\tau \intO{ \left[ \left< Y_{t,x}; [\vr]_{\rm ess} \right> s(r, \Theta) \partial_t \Theta + \left< Y_{t,x}; [\vc{m}]_{\rm ess} \right> \cdot s(r, \Theta) \Grad \Theta    \right] } \dt\\
- & \int_0^\tau \intO{ \left[ \left< Y_{t,x} ; \left[ \vr Z\left( s(\vr,E) \right) \right]_{\rm res} \right>  \partial_t \Theta +
\left< Y_{t,x} ; \left[ Z \left( s(\vr,E) \right) \vc{m} \right]_{\rm res} \right> \cdot \Grad \Theta  \right] } \dt\\
+ &\int_0^\tau \intO{ \left[ \left< Y_{t,x}; [\vr]_{\rm res} \right> s(r, \Theta) \partial_t \Theta + \left< Y_{t,x}; [\vc{m}]_{\rm res} \right> \cdot s(r, \Theta) \Grad \Theta    \right] } \dt
\end{split}
\]
\[
\begin{split}
=& \int_0^\tau \intO{ \left< Y_{t,x}; [ \vr (s(r, \Theta) - s(\vr,E) )]_{\rm ess}  \right> \partial_t \Theta } \ \dt \\
+& \int_0^\tau \intO{ \left< Y_{t,x}; [ \vc{m} (s(r, \Theta) - s(\vr, E))]_{\rm ess}  \right> \cdot \Grad \Theta } \ \dt
\\
- & \int_0^\tau \intO{ \left[ \left< Y_{t,x} ; \left[ \vr Z\left( s(\vr,E) \right) \right]_{\rm res} \right>  \partial_t \Theta +
\left< Y_{t,x} ; \left[ Z \left( s(\vr,E) \right) \vc{m} \right]_{\rm res} \right> \cdot \Grad \Theta  \right] } \dt\\
+ &\int_0^\tau \intO{ \left[ \left< Y_{t,x}; [\vr]_{\rm res} \right> s(r, \Theta) \partial_t \Theta + \left< Y_{t,x}; [\vc{m}]_{\rm res} \right> \cdot s(r, \Theta) \Grad \Theta    \right] } \dt,
\end{split}
\]
 where the residual terms are controlled in view of (\ref{WS8}).

As for the essential components, we may pass to the original variables $(\vr, \vt)$ to observe that
\[
[s(\vr, \vt(\vr,E)) - s(r, \Theta)]_{\rm ess} \approx \frac{\partial s(r, \Theta) }{\partial \vr} [\vr - r]_{\rm ess} +
\frac{\partial s(r, \Theta) }{\partial \vt} [\vt(\vr,E) - \Theta]_{\rm ess},
\]
where the difference proportional to 
\[
[ \vr - r ]^2_{\rm ess} + [E - r e(r, \Theta) ]^2_{\rm ess},
\]
is absorbed by the left-hand side of (\ref{WS9}).

Summing up the previous discussion, we may replace (\ref{WS9}) by
\begin{equation} \label{WS11}
\begin{split}
& \intO{ \left< Y_{\tau,x} ; \mathcal{E} \left(\vr, E, \vc{m} \ \Big| r, \Theta, \vc{U} \right) \right> }
+ \mathcal{D}(\tau)
\\
&\aleq - \int_0^\tau \intO{ \left< Y_{t,x}; \vr \left( \frac{\partial s(r, \Theta) }{\partial \vr} [\vr - r]_{\rm ess} +
\frac{\partial s(r, \Theta) }{\partial \vt} [\vt(\vr, E) - \Theta]_{\rm ess} \right) \right> \partial_t \Theta } \ \dt
\\
&- \int_0^\tau \intO{ \left< Y_{t,x}; \vc{m} \left( \frac{\partial s(r, \Theta) }{\partial \vr} [\vr - r]_{\rm ess} +
\frac{\partial s(r, \Theta) }{\partial \vt} [\vt(\vr,E) - \Theta]_{\rm ess} \right) \right> \cdot \Grad \Theta } \ \dt
\\
&+ \int_0^\tau \intO{ \left[p(r, \Theta) \Div \vc{U} - \left< Y_{t,x}; p(\vr,E) \right> \Div \vc{U}  \right] } \dt \\
&+ \int_0^\tau \intO{ \left[ \left< Y_{t,x} ; r - \vr \right> \frac{1}{r} \partial_t p(r, \Theta) -
\left< Y_{t,x} ; \vr \vc{U} \right> \cdot \frac{1}{r} \Grad p(r, \Theta) + \Grad p(r, \Theta) \cdot \vc{U}   \right] } \ \dt\\
& + \int_0^\tau \left[ \intO{ \left< Y_{t,x} ; \mathcal{E} \left(\vr, E, \vc{m} \ \Big| r, \Theta, \vc{U} \right) \right> }
+ \mathcal{D}(t) \right] \dt \
\mbox{for a.a.}\ \tau \in (0,T).
\end{split}
\end{equation}

\medskip

{\bf Step 3}

Using the fact that $r$, $\vc{U}$ satisfy the equation of continuity we get, after a tedious but straightforward manipulation, the identity
\[
\begin{split}
&(r - \vr) \frac{1}{r} \partial_t p(r, \Theta) + \Grad p(r, \Theta) \cdot \vc{U} - \frac{\vr}{r} \vc{U} \cdot \Grad p(r, \Theta) + \Div \vc{U} (p(r,\Theta) - p(\vr, \vt) )\\
&= \Div \vc{U} \left( p(r, \Theta) - \frac{\partial p(r,\Theta)}{\partial \vr} (r - \vr) - \frac{\partial p(r,\Theta)}{\partial \vt} (\Theta - \vt)
- p(\vr, \vt)
\right)\\
& + r (\vr - r) \frac{\partial s(r, \Theta) }{\partial \vr} \left( \partial_t \Theta + \vc{U} \cdot \Grad \Theta \right) +
r (\vt - \Theta) \frac{\partial s(r, \Theta) }{\partial \vt} \left( \partial_t \Theta + \vc{U} \cdot \Grad \Theta \right).
\end{split}
\]
In view of hypothesis (\ref{WS1a}) the residual part of the expression on the left-hand side is controlled and we may go back to (\ref{WS11}) to deduce the
desired conclusion: 
\begin{equation} \label{WS12}
\begin{split}
& \intO{ \left< Y_{\tau,x} ; \mathcal{E} \left(\vr, E, \vc{m} \ \Big| r, \Theta, \vc{U} \right) \right> }
+ \mathcal{D}(\tau)
\\
&\aleq
\int_0^\tau \intO{ \left< Y_{t,x}; \left[p(r, \Theta) - \frac{\partial p(r,\Theta)}{\partial \vr} (r - \vr) - \frac{\partial p(r,\Theta)}{\partial \vt} (\Theta - \vt(\vr,E))
- p(\vr, \vt(\vr,E)) \right]_{\rm ess} \right> \Div \vc{U} } \ \dt  \\
& + \int_0^\tau \left[ \intO{ \left< Y_{t,x} ; \mathcal{E} \left(\vr, E, \vc{m} \ \Big| r, \Theta, \vc{U} \right) \right> }
+ \mathcal{D}(t) \right] \dt
\\
&\aleq
\int_0^\tau \intO{ \left< Y_{t,x}; \left[|\vr - r |^2 + |E -re(r,\Theta)|^2\right]_{\rm ess} \right>  } \ \dt  \\
& + \int_0^\tau \left[ \intO{ \left< Y_{t,x} ; \mathcal{E} \left(\vr, E, \vc{m} \ \Big| r, \Theta, \vc{U} \right) \right> }
+ \mathcal{D}(t) \right] \dt \
\mbox{for a.a.}\ \tau \in (0,T).
\end{split}
\end{equation}

Applying Gronwall' lemma we deduce that the left-hand side of (\ref{WS12}) vanishes for a.a. $\tau \in (0,T)$.

We have shown the following result.

\Cbox{Cgrey}{

\begin{Theorem} {\bf[Weak (measure-valued) - strong uniqueness principle]} \label{TWS1}

Let the thermodynamic functions $e = e(\vr, \vt)$, $s = s(\vr, \vt)$, and $p = p(\vr, \vt)$ satisfy Gibbs' relation (\ref{i5}), the hypothesis
of thermodynamic stability (\ref{WS1}), and let
\begin{equation} \label{HYP}
|p(\vr, \vt)| \leq c (1 + \vr + \vr|s(\vr, \vt)| + \vr e(\vr, \vt) ).
\end{equation}

Let $[r,\Theta, \vc{U}]$ be a continuously differentiable classical solution of the Euler system (\ref{i1}--\ref{i3}) in $(0,T) \times \Omega$ starting
from the initial data $(r_0, \Theta_0, \vc{U}_0)$ satisfying
\[
r_0 > 0, \ \Theta_0 > 0.
\]
Assume that $[ Y_{t,x}; \mathcal{D} ]$ is a dissipative measure valued solution of the same problem in the sense specified in Definition \ref{DMV1}
such that
\[
Y_{0,x} = \delta_{[ r_0(x), r_0 e(r_0, \Theta_0)(x), r_0 \vc{U}_0 (x) ]} \ \mbox{for a.a.}\ x \in \Omega.
\]

Then $\mathcal{D} = 0$ and
\[
Y_{t,x} = \delta_{[ r(t, x), r e(r, \Theta)(t,x), r \vc{U} (t, x) ]} \ \mbox{for any}\ (t,x) \in (0,T) \times \Omega.
\]

\end{Theorem}

}

\section{Conclusion}
\label{C}

We have introduced the concept of \emph{dissipative measure-valued solution} to the complete Euler system (\ref{i1}--\ref{i4}). Such a solution appears as a natural cluster
point of families of weak solutions or their viscous approximations. We expect also certain numerical schemes to generate this kind of solutions, cf.
\cite{FeiLuk}. The main result stated in Theorem \ref{TWS1} above asserts that a dissipative measure-valued solution coincides with a strong solution starting from the same initial data on the life span of the latter. In particular, if this is the case, any sequence generating the measure-valued solution necessarily converges pointwise to the strong solution.
Such a result can be used for proving convergence of certain numerical schemes as in the simpler barotropic case discussed in \cite{FeiLuk}.

The fact that the measure--valued formulation (\ref{MV22}--\ref{MV27}) contains the cut-off function $Z$ may seem restrictive although quite natural in the present context.
The cut off can be dropped, meaning taking $Z(s) = s$, provided sufficiently strong {\it a priori} bounds are available to control integrability of $\vr s \vu = s \vc{m}$.
As we have seen in Section \ref{MV}, these bounds followed from boundedness from below of the entropy of the system. Similar bounds can be obtained directly from the
energy balance provided the constitutive thermodynamic functions satisfy a technical restriction
\begin{equation} \label{C1}
\vr |s(\vr, \vt)|^2 \aleq (1 + \vr + \vr e(\vr, \vt) ).
\end{equation}
It can be shown that (\ref{C1}) holds for a general  monoatomic gas satisfying the caloric equation of state
\[
p = \frac{2}{3} \vr e,
\]
provided that the associated entropy $s = s(\vr, \vt)$ complies with the Third law of thermodynamics, specifically,
\[
\lim_{\vt \to 0} s(\vr, \vt) = 0 \ \mbox{for any} \ \vr > 0,
\]
cf. \cite[Chapter 1, Part 1.4]{FENO6}.

\def\cprime{$'$} \def\ocirc#1{\ifmmode\setbox0=\hbox{$#1$}\dimen0=\ht0
  \advance\dimen0 by1pt\rlap{\hbox to\wd0{\hss\raise\dimen0
  \hbox{\hskip.2em$\scriptscriptstyle\circ$}\hss}}#1\else {\accent"17 #1}\fi}

%\bibliographystyle{plain}
%\bibliography{citace}

\begin{thebibliography}{10}

\bibitem{BALL2}
J.M. Ball.
\newblock A version of the fundamental theorem for {Y}oung measures.
\newblock {\em In Lect. Notes in Physics 344, Springer-Verlag}, pages 207--215,
  1989.

\bibitem{BAMU}
J.M. Ball and F.~Murat.
\newblock Remarks on {C}hacons biting lemma.
\newblock {\em Proc. Amer. Math. Soc.}, {\bf 107}:655--663, 1989.

\bibitem{BenSer}
S.~Benzoni-Gavage and D.~Serre.
\newblock {\em Multidimensional hyperbolic partial differential equations,
  {F}irst order systems and applications}.
\newblock Oxford Mathematical Monographs. The Clarendon Press Oxford University
  Press, Oxford, 2007.

\bibitem{BEROFO}
S.E.~Bechtel, F.J.~Rooney, and M.G.Forest.
\newblock Connection between stability, convexity of internal energy,
and the second law for compressible {N}ewtonian fuids.
\newblock{\em J. Appl. Mech.}, {\bf 72}:299-300, 2005.


\bibitem{BrDeSz}
Y.~Brenier, C.~De~Lellis, and L.~Sz{\'e}kelyhidi, Jr.
\newblock Weak-strong uniqueness for measure-valued solutions.
\newblock {\em Comm. Math. Phys.}, {\bf 305}(2):351--361, 2011.


\bibitem{CheFr2}
G.-Q.~Chen and H.~Frid.
\newblock Uniqueness and asymptotic stability of {R}iemann solutions for
              the compressible {E}uler equations.
\newblock{\em Trans. Amer. Math. Soc.}, {{\bf 353}} (3):1103--1117, 2001.


\bibitem{Chiod}
E.~Chiodaroli.
\newblock A counterexample to well-posedness of entropy solutions to the
  compressible {E}uler system.
\newblock {\em J. Hyperbolic Differ. Equ.}, {\bf 11}(3):493--519, 2014.

\bibitem{ChiDelKre}
E.~Chiodaroli, C.~{D}e {L}ellis, and O.~Kreml.
\newblock Global ill-posedness of the isentropic system of gas dynamics.
\newblock {\em Comm. Pure Appl. Math.}, {\bf 68}(7):1157--1190, 2015.

\bibitem{DelSze3}
C.~De~Lellis and L.~Sz{\'e}kelyhidi, Jr.
\newblock On admissibility criteria for weak solutions of the {E}uler
  equations.
\newblock {\em Arch. Ration. Mech. Anal.}, {\bf 195}(1):225--260, 2010.

\bibitem{DiP2}
R.J. Di{P}erna.
\newblock Measure-valued solutions to conservation laws.
\newblock {\em Arch. Rat. Mech. Anal.}, {\bf 88}:223--270, 1985.

\bibitem{FGSWW1}
E.~Feireisl, P.~Gwiazda, A.~{\' S}wierczewska-{G}wiazda, and E.~Wiedemann.
\newblock Dissipative measure-valued solutions to the compressible
  {N}avier-{S}tokes system.
\newblock {\em Calc. Var. Partial Differential Equations} {\bf 55}(6), Art. 141, 20pp., 2016. 




\bibitem{FeiLuk}
E.~Feireisl and M.~Luk{\' a}{\v c}ov{\' a}-{M}edvi{\v d}ov{\' a}.
\newblock Convergence of a mixed finite element finite volume scheme for the
  isentropic {N}avier--{S}tokes system via dissipative measure-valued
  solutions.
\newblock 2016.
\newblock Peprint no. IM 2016-43 IM Prague.

\bibitem{FENO6}
E.~Feireisl and A.~Novotn{\' y}.
\newblock {\em Singular limits in thermodynamics of viscous fluids}.
\newblock Birkh{\" a}user-Verlag, Basel, 2009.

\bibitem{FeiNov10}
E.~Feireisl and A.~Novotn{\' y}.
\newblock Weak-strong uniqueness property for the full
  {N}avier-{S}tokes-{F}ourier system.
\newblock {\em Arch. Rational Mech. Anal.}, {\bf 204}:683--706, 2012.

\bibitem{FjKaMiTa}
U.~K. Fjordholm, R.~K{\" a}ppeli, S.~Mishra, and E.~Tadmor.
\newblock Construction of approximate entropy measure valued solutions for
  hyperbolic systems of conservation laws.
\newblock {\em Foundations Comp. Math.}, pages 1--65, 2015.

\bibitem{FjMiTa2}
U.~S. Fjordholm, S.~Mishra, and E.~Tadmor.
\newblock Arbitrarily high-order accurate entropy stable essentially
  nonoscillatory schemes for systems of conservation laws.
\newblock {\em SIAM J. Numer. Anal.}, {\bf 50}(2):544--573, 2012.

\bibitem{FjMiTa1}
U.~S. Fjordholm, S.~Mishra, and E.~Tadmor.
\newblock On the computation of measure-valued solutions.
\newblock {\em Acta Numer.}, {\bf 25}:567--679, 2016.

\bibitem{GSWW}
P.~Gwiazda, A.~{\' S}wierczewska-{G}wiazda, and E.~Wiedemann.
\newblock Weak-strong uniqueness for measure-valued solutions of some
  compressible fluid models.
\newblock  {\em Nonlinearity}, {\bf 28}(11):3873–-3890, 2015.

\bibitem{KrZa}
D.~Kr{\"o}ner and W.~M. Zajaczkowski.
\newblock Measure-valued solutions of the {E}uler equations for ideal
  compressible polytropic fluids.
\newblock {\em Math. Methods Appl. Sci.}, {\bf 19}(3):235--252, 1996.

\bibitem{MNRR}
J.~M{\' a}lek, J.~Ne{\v c}as, M.~Rokyta, and M.~R{\accent23u}{\v z}i{\v c}ka.
\newblock {\em Weak and measure-valued solutions to evolutionary PDE's}.
\newblock Chapman and Hall, London, 1996.

\end{thebibliography}

\end{document}